\documentclass[10pt]{studiamnew}
\usepackage{graphicx}
\usepackage{amsmath}
\usepackage{amssymb}
\newtheorem{theorem}{Theorem}[section]
\newtheorem{lemma}[theorem]{Lemma}

\theoremstyle{definition}
\newtheorem{remark}[theorem]{Remark}

\renewcommand{\theequation}{\thesection.\arabic{equation}}
\numberwithin{equation}{section}

\begin{document}
%
%

\title[Approximation   theorems for multivariate Taylor-Abel-Poisson means]
{Approximation   theorems  for  multivariate Taylor-Abel-Poisson means}
\author{J\"{u}rgen Prestin}
\address{University of L{\"{u}}beck \\ Institute of Mathematics\\
 Ratzeburger Allee 160\\
23562 L\"{u}beck \\
Germany}
\email{prestin@math.uni-luebeck.de}
\author{Viktor Savchuk}
\address{Institute of Mathematics of the National Academy of Sciences of Ukraine\\
3, Tereschenkivska street\\
01004  Kyiv\\
Ukraine}
\email{vicsavchuk@gmail.com}
\author{Andrii Shidlich}
\address{Institute of Mathematics of the National Academy of Sciences of Ukraine\\
3, Tereschenkivska street\\
01004  Kyiv\\
Ukraine}
\email{shidlich@gmail.com}

%
\subjclass{41A27, 42A45, 	41A35}
\keywords{direct approximation  theorem;  inverse approximation  theorem; Taylor-Abel-Poisson means, $K$-functional, multiplier}
\begin{abstract}
We obtain direct and inverse  approximation  theorems of  functions of several variables by Taylor-Abel-Poisson means in the integral metrics. We also show that norms of multipliers in the spaces $L_{p,Y}(\mathbb T^d)$ are equivalent  for all positive integers $d.$
\end{abstract}
\maketitle

\section{Introduction} It is well-known that any function $f \in L_p(\mathbb T^1)$  that is different from a constant can be approximated by its Abel-Poisson means $f(\varrho,\cdot)$ with a precision not better than $1-\varrho$. It relates to the so-called saturation property of this approximation method. From this property, it follows that for any $f\in L_p(\mathbb T^1)$, the relation $\|f-f(\varrho,\cdot)\|_{_{\scriptstyle p}} =\mbox{\tiny  $\mathcal O$}(1-\varrho)$, $\varrho\to 1-$, only holds  in the trivial case when  $f$ is a constant  function. Therefore, any additional restrictions on the smoothness of functions do not give us any order of approximation better than  $1-\varrho$. In this connection, a natural question is to find a linear operator,  constructed similarly to the Poisson operator, which takes into account the smoothness properties of functions and at the same time, for a given functional class,  is the best in a certain sense.
In [\ref{Savchuk_Zastavnyi}], for classes of convolutions whose kernels were generated by some moment sequences, the authors proposed a general method of construction of similar operators that take into account properties of such kernels and hence, the smoothness of functions from corresponding  classes. One example of such operators are the operators $ A_{\varrho, r} $, which are the main subject of study in this paper.

The operators $A_{\varrho,r}$ were first studied in [\ref{Savchuk}] where, in the terms of these operators, the author gave the structural characteristic of Hardy-Lipschitz classes $H^r_p\mathop{\rm Lip}\alpha$  of
functions of one variable,
holomorphic
on
the unit disc
of
the complex plane. In  [\ref{Savchuk_Shidlich}],  in
terms of approximation estimates of such operators in
some spaces $S^p$ of Sobolev type, the authors give a constructive description of classes of functions  of several variables whose generalized derivatives belong to the classes $S^pH_\omega$.  In [\ref{Prestin_Savchuk_Shidlich}], direct and inverse approximation theorems of $2\pi$-periodic functions by the operators $A_{\varrho, r}$ were given in the terms of $K$--functionals of functions generated by their radial derivatives.

Approximations of functions of one variable by similar operators of polynomial type  were studied in [\ref{Leis}], [\ref{Butzer_Sunouchi}], [\ref{Chui}], [\ref{Holland}], [\ref{Mohapatra}], [\ref{Chandra}] etc. In particular, in [\ref{Chui}], the authors found the  degree of convergence of the well-known Euler and Taylor means to the functions $f$ from some subclasses of the Lipschitz classes ${\rm Lip} \alpha$ in the uniform norm. In [\ref{Mohapatra}], the analogous results for Taylor means were obtained in the $L_p$--norm.

In the present paper, we continue the study of  approximative properties of  the operators $A_{\varrho, r}$. In particular, we extend the results of the paper  [\ref{Prestin_Savchuk_Shidlich}] to the  multivariate case and prove direct and inverse  approximation  theorems of  functions of several variables by the operators $A_{\varrho,r}$  in the integral metrics.  We also show that norms of multipliers in the spaces $L_{p,Y}(\mathbb T^d)$ are equivalent for all positive integers $d.$

\section{Preliminaries} Let $d$ be an integer, let $\mathbb R^d$, $\mathbb R^d_+$ and $\mathbb Z^d$ be the sets of all vectors ${\bf k}: = (k_1, \ldots, k_d)$ with  real,   real non-negative  and   integer  coordinates  respectively. Set  $\mathbb T^d:={\mathbb R}^d/ 2\pi {\mathbb Z}^d$.

Further, let $L_p(\mathbb T^d)$, $1\le p\le \infty,$  be the space of all functions  $f({\bf x})=f(x_1,\ldots,x_d)$ defined on  $\mathbb R^d$,  $2\pi$-periodic in each variable with the finite norm
\begin{equation}\label{norm Lp}
\|f\|_{_{\scriptstyle p}}=\|f\|_{_{\scriptstyle L_p(\mathbb T^d) }}:=\left\{\begin{matrix}\Big(\displaystyle{\int_{\mathbb T^d}}|f({\bf x})|^p{\mathrm d}\sigma({\bf x})\Big)^{\frac 1p},\quad \hfill & 1\le p< \infty,\cr
\mathop{\rm ess\,sup}_{{\bf x}\in \mathbb T^d}|f({\bf x})|,\quad \hfill &  p=\infty,\end{matrix}\right.
\end{equation}
where $\sigma$ is the normalized Lebesgue measure on  $\mathbb T^d$.

Let $({\bf x},{\bf y}):=x_1y_1+\ldots+x_dy_d$ denote the inner product of the elements ${\bf x},{\bf y}\in {\mathbb R}^d$. Let us set ${\mathrm e}_{\bf k}:={\mathrm e}_{\bf k}({\bf x})={\mathrm e}^{{\mathrm i}({\bf k},{\bf x})}$, ${\bf k}\in {\mathbb Z}^d$, and for any function  $f\in L_1(\mathbb T^d)$, define its Fourier coefficients by
$$
\widehat f_{\bf k}:=\int_{\mathbb T^d}f({\bf x})\overline{\mathrm e}_{\bf k}({\bf x}){\mathrm d}\sigma({\bf x}),\quad{\bf k}\in\mathbb Z^d,
$$
where $\overline{z}$ is the complex-conjugate number of $z$.

Set $|{\bf k}|_1:=\sum_{j=1}^d |k_j|$,  and  for any function $f\in L_1(\mathbb T^d)$ with the Fourier series of the form
 \begin{equation}\label{Fourier series}
S[f]({\bf x})=
\sum_{{\bf k}\in {\mathbb Z}^d}\widehat f_{\bf k} {\mathrm e}_{\bf k}({\bf x})
 =\sum_{\nu=0}^\infty\sum_{|{\bf k}|_1=\nu}\widehat f_{\bf k} {\mathrm e}_{\bf k}({\bf x}),
\end{equation}
denote by $f\left(\mbox{\boldmath${\bf \varrho}$},{\bf x}\right)$ its Poisson integral  (the Poisson operator), i.e.,
 \begin{equation}\label{Poisson operator}
f\left(\mbox{\boldmath${\bf \varrho}$},{\bf x}\right):=\int_{\mathbb T^d}
f({\bf x}+{\bf s})P(\mbox{\boldmath${\bf \varrho}$},{\bf s}){\mathrm d}\sigma({\bf s}),
\end{equation}
where  $\mbox{\boldmath${\bf \varrho}$}\in \mathbb R^d_+$, ${\bf x}\in {\mathbb R}^d$, the function $P(\mbox{\boldmath${\bf \varrho}$},{\bf x}):=\sum_{{\bf k}\in {\mathbb Z}^d}
\mbox{\boldmath${\bf \varrho}$}^{|{\bf k}|} {\mathrm e}_{\bf k}({\bf x})$ is the  Poisson kernel,
$\mbox{\boldmath${\bf \varrho}$}^{|{\bf k}|}:=\varrho_1^{|k_1|}\cdot\cdot\cdot \varrho_d^{|k_d|}$.

In what follows,  the expression $f(\varrho,{\bf x})$ means the Poisson integral, where {\boldmath${\bf \varrho}$} is a vector with the same coordinates, i.e.,  {\boldmath${\bf \varrho}$}$\,=(\varrho,\ldots,\varrho).$ In such case, we have $P(\varrho,{\bf x}):=\sum_{\nu=0}^{\infty} \varrho^{\nu}\sum_{|{\bf k}|_1=\nu}  {\mathrm e}_{\bf k}({\bf x})$.

Let $f\in L_1(\mathbb T^d)$. For  $\varrho\in [0,1)$ and $r\in\mathbb N$, we set
\begin{equation}\label{def Ar}
A_{\varrho,r}(f)({\bf x}):=
 \sum_{\nu=0}^{\infty}\lambda_{\nu,r}(\varrho)\sum_{|{\bf k}|_1=\nu}\widehat f_{\bf k}{\mathrm e}_{\bf k}({\bf x}),
\end{equation}
where for $\nu=0,1,\ldots,r-1$,   the numbers $\lambda_{\nu,r}(\varrho)\equiv 1$ and for
$\nu=r, r+1,\ldots$,
\begin{equation}\label{lambda for H^r}
 \lambda_{\nu,r}(\varrho):=\sum_{j=0}^{r-1}
{\nu\choose j}(1-\varrho)^{j}\varrho^{\nu-j}
=
\sum_{j=0}^{r-1}\frac{(1-\varrho)^{j}}{j!}~\frac{d^j}{d\varrho^j}\varrho^{\nu}.
\end{equation}
The transformation  $A_{\varrho,r}$ can be considered as a linear operator
on  $L_1(\mathbb T^d)$ into itself. Indeed, $\lambda_{\nu,r}(0)=0$ and for all  $\nu=r,r+1,\ldots$ and $\varrho\in(0,1)$, 
\[
\sum_{j=0}^{r-1}
{\nu\choose j}(1-\varrho)^j\varrho^{\nu-j}\le
rq^{\nu}\nu^{r-1},~\mbox{where}~0<q:=\max\{1-\varrho,\varrho\}<1.
\]
Therefore, for any function $f\in L_1(\mathbb T^d)$ and for any $0<\varrho<1$, the series on the right-hand side of (\ref{def Ar}) is majorized by the convergent series $2r\|f\|_{_{\scriptstyle 1}}\sum_{\nu=r}^{\infty}q^{\nu}\nu^{r-1}$.


Leis [\ref{Leis}] considered for  $f\in L_p(\mathbb T^1)$, $1<p<\infty$, the transformation
\[
L_{\varrho, r}(f)(x):=\sum_{k=0}^{r-1}\frac{ \mathrm{d}^k f(x)}{ \mathrm{d} n^k}\cdot\frac{(1-\varrho)^k}{k!},\quad r\in\mathbb N,
\]
where
\[
\frac{ \mathrm{d} f(x)}{ \mathrm{d} n}=\left.-\frac{\partial f(\varrho,x)}{\partial\varrho}\right|_{\varrho=1}
\]
is the normal derivative of the function $f$. He showed that if $1<p<\infty$ and
\[
\|f(\varrho, \cdot)-L_{\varrho, r}(f)(\cdot)\|_{_{\scriptstyle p}} =
{\mathcal O}\Big({(1-\varrho)^r}/{r!}\Big),\quad\varrho\to 1-,
\]
then  $ \mathrm{d}^r f/ \mathrm{d} n^r \in L_p(\mathbb T^1)$.

 Butzer and Sunouchi [\ref{Butzer_Sunouchi}] considered for $f\in L_p(\mathbb T^1)$, $1\le p<\infty$, the transformation
\[
B_{\varrho,r}(f)(x):=\sum_{k=0}^{r-1}(-1)^{\frac{k+1}{2}}f^{\{k\}}(x)\frac{(-\ln\varrho)^k}{k!},
\]
where $f^{\{k\}}:=f^{(k)}$ for $k\in 2\mathbb Z_+$ and $f^{\{k\}}:=\widetilde f^{(k)}$ for $k-1
\in 2\mathbb Z_+$, where
 $$
 \widetilde f(x)=\lim\limits_{\varepsilon\to 0}-\frac 1\pi\int\limits_\varepsilon^\pi (f(x+u)-f(x-u))\frac 12 {\cot}\frac u2 {\mathrm d}u.
 $$
They proved the following theorem:

{\bf Theorem A [\ref{Butzer_Sunouchi}].} {\it Assume that  $f\in L_p(\mathbb T^1)$, $1\le p<\infty$.

$i)$
If
the derivatives $f^{\{j\}},$ $j=0,1,\ldots, r-1,$ are absolutely continuous  and $f^{\{r\}}\in L_p(\mathbb T^1)$, then
\begin{equation}\label{B-S approx}
\|f(\varrho,\cdot)-B_{\varrho, r}(f)(\cdot)\|_{_{\scriptstyle p}}={\mathcal O}\Big({(-\ln\varrho)^r}/{r!}\Big),\quad\varrho\to 1-~%
.
\end{equation}

$ii)$
If
the derivatives $f^{\{j\}},$ $j=0,1,\ldots, r-2,$ $r\ge 2,$  are absolutely continuous, $f^{\{r-1\}}\in L_p(\mathbb T^1)$, $1<p<\infty,$ and relation (\ref{B-S approx}) holds, then $\widetilde f^{\{r-1\}}$ is absolutely continuous and $\widetilde f^{\{r\}}\in L_p(\mathbb T^1)$.
}

These results
summarize the approximation behaviour
of the operators $L_{\varrho, r}$ and $B_{\varrho, r}$ in the space $L_p(\mathbb T^1)$. 
In particular, Leis's result and the statement  $ii)$ of Theorem A
represent the so-called 
inverse theorems and the statement $i)$ is the
so-called 
direct theorem.
Direct and inverse theorems are among the main  theorems of approximation theory. They were studied by many authors. Here, we mention only the books [\ref{Butzer_Nessel}, \ref{DeVore Lorentz}, \ref{Trigub_Bellinsky}] which contain fundamental results in this subject.
The result of Leis and Theorem A are based on the investigations
in
the papers [\ref{Butzer_Tillmann}, \ref{Butzer}], where the authors find the direct and inverse approximation theorems for the one-parameter  semi-groups of bounded
linear transformations  $\{T(t)\}$ of
some 
Banach space $X$ into
itself by  the ``Taylor polynomial''  
$\sum_{k=0}^{r-1}(t^k/k!)A^kf$, where $Af$ is the infinitesimal operator of a semi--group $\{T(t)\}$.

The transformations $A_{\varrho, r}$ considered in
this
paper are similar to the transformations $L_{\varrho, r}$ and $B_{\varrho,r}$ as they are also  based on the ``Taylor polynomials''. 
The
relation between the operators $A_{\varrho, r}$ and  the ``Taylor polynomials''  is
 shown in the following statement.


{\lemma\label{Lemma1}   Assume that  $f\in L_1(\mathbb T^d)$. Then for any numbers $r\in \mathbb N,$ $\varrho\in[0,1)$ and ${\bf x}\in  \mathbb T^d$,
\begin{equation}\label{A_P}
A_{\varrho, r}(f)({\bf x})=\sum_{j=0}^{r-1}\frac{\partial^j f\left(\varrho,{\bf x}\right)}{\partial\varrho^j}\cdot\frac{(1-\varrho)^j}{j!}.
\end{equation}
}

\noindent {\it Proof.} With respect to the variable $\varrho$, let us differentiate  the  decomposition of the Poisson integral into the uniformly convergent series
\begin{equation}\label{series fo Poisson}
f\left(\varrho,{\bf x}\right)=\sum_{\nu=0}^{\infty}\varrho^{\nu}\sum_{|{\bf k}|_1=\nu}\widehat f_{\bf k}{\mathrm e}_{\bf k}({\bf x}), \quad    \varrho\in[0,1),~
{\bf x}\in  \mathbb T^d.
\end{equation}
We see that for any $j=0,1,\ldots$
\begin{equation}\label{Poisson_Derivative}
 \frac{\partial^j f\left(\varrho,{\bf x}\right)}{\partial\varrho^j}= \sum_{\nu=j}^{\infty}\frac{\nu!}{(\nu-j)!}
 \varrho^{\nu-j}\sum_{|{\bf k}|_1=\nu}\widehat f_{\bf k}{\mathrm e}_{\bf k}({\bf x}) .
\end{equation}
Since $\sum_{j=0}^{\nu}
{\nu\choose j}(1-\varrho)^j\varrho^{\nu-j}=\big((1-\varrho)+\varrho\big)^{\nu}=1,~\nu=0,1,\ldots,
$
then
\[
\sum_{j=0}^{r-1}\frac{\partial^j f\left(\varrho,{\bf x}\right)}{\partial\varrho^j}
\cdot\frac{(1-\varrho)^j}{j!}=
\sum_{\nu=0}^{r-1}\sum_{j=0}^{\nu}{\nu\choose j}(1-\varrho)^j\varrho^{\nu-j} \sum_{|{\bf k}|_1=\nu}\widehat f_{\bf k}{\mathrm e}_{\bf k}({\bf x})
\]
\[
+\sum_{\nu=r}^{\infty}\sum_{j=0}^{r-1}
{\nu\choose j}(1-\varrho)^{j}\varrho^{\nu-j}\sum_{|{\bf k}|_1=\nu}\widehat f_{\bf k}{\mathrm e}_{\bf k}({\bf x})=A_{\varrho, r}(f)({\bf x}).
\]
\vskip -3mm$\hfill\Box$

\section{Direct and inverse  approximation  theorems}

\subsection{Radial derivatives and $K$-functionals}\label{Subsetion2.1}

If for a function $f\in L_1(\mathbb T^d)$ and for a positive integer  $n$ there exists the function $g\in L_1(\mathbb T^d)$  such that
 \[
\widehat g_{\bf k}=\left\{\begin{matrix} 0,\hfill  & \quad \ \mbox{\rm if}\quad |{\bf k}|_1=\nu<n,\quad \  \cr  {\displaystyle \frac{\nu!}{(\nu-n)!}}\widehat f_{\bf k},\quad   \hfill & \mbox{\rm if}\quad |{\bf k}|_1=\nu \ge n,\end{matrix}\right.\quad  {\bf k}\in {\mathbb Z}^d,\quad \nu=0,1,\ldots,
\]
then we say that for the function $f$, there exists the radial derivative  $g$ of order $n$ for which we use the notation $f^{[n]}$.

Let us note that if the function $f^{[r]}\in L_1(\mathbb T^d)$, then its Poisson integral can be presented
as
\begin{equation}\label{diff f[n]}
f^{[r]}(\varrho,{\bf x})=(f(\varrho,\cdot))^{[r]}({\bf x})=\varrho^r\frac{\partial^r f\left(\varrho,{\bf x}\right)}{\partial\varrho^r}
\quad \varrho\in[0,1),~\forall~{\bf x}\in {\mathbb T}^d.
\end{equation}


In the space $L_p(\mathbb T^d)$,
the 
$K$--functional of
a
function $f$ (see, for example, [\ref{DeVore Lorentz}, Chap.~6])  generated by the radial derivative of order  $n$  is the following quantity:
\[
K_n(\delta, f)_p:=\inf\Big\{\left\|f-h\right\|_p+\delta^n \|h^{[n]} \|_p: h^{[n]}\in L_p(\mathbb T^d)\Big\},\quad\delta>0.
\]

\subsection{Main results}\label{Subsetion2.2}


Let  $\mathbb Z^d_-$  denote the set of all vectors ${\bf k}: = (k_1, \ldots, k_d)$ with negative integer coordinates, $\mathbb Z^d_+:=\mathbb Z^d\cap \mathbb R^d_+$ and  $Y:={\mathbb Z}^d_+\cup {\mathbb Z}^d_-$. Let also $L_{p,Y}(\mathbb T^d)$ be the set of all functions $f$ from $L_p(\mathbb T^d)$ such that the Fourier coefficients $\widehat f_{\bf k}=0$ for all ${\bf k}\in {\mathbb Z}^d\setminus Y$.  Further, we consider the functions $\omega(t)$, $t\in [0,1]$, satisfying the following conditions 1)--\,4):
\mbox{{\bf 1)} $\omega(t)$ }is continuous on $[0,1]$; {\bf 2)} $\omega(t)$ is monotonically increasing; {\bf 3)}~$\omega(t){\not=}\,0$ for all $t\in (0,1]$; {\bf 4)} $\omega(t)\to 0$ as $t\to 0$; and the  well-known Zygmund--Bari--Stechkin conditions (see, for example, [\ref{Bari_Stechkin}]):
$$
({\mathcal Z}): \int_0^\delta\frac{\omega(t)}{t}{\mathrm d}t={\mathcal O}(\omega(\delta)),\ \
({\mathcal Z}_n):   \int_\delta^1\frac{\omega(t)}{t^{n+1}}{\mathrm d}t={\mathcal O}\Big(\frac{\omega(\delta)}{\delta^n}\Big),\ n\in {\mathbb N},\, \delta\to 0+.
$$


The main results of this paper are contained in the following two statements:

{\theorem\label{Th1} Assume that $f\in L_{p,Y}(\mathbb T^d),$ $1\le p\le\infty$, $n, r\in\mathbb N$, $n\le r$ and the function  $\omega(t)$, $t\in [0,1]$, satisfies conditions 1)--4) and $({\mathcal Z})$ . If
\begin{equation}\label{K-funct est}
f^{[r-n]}\in L_p(\mathbb T^d)\quad and\quad K_{n}\left(\delta, f^{[r-n]}\right)_p={\mathcal O}(\omega(\delta)),\ \delta\to 0+,
\end{equation}
then
\begin{equation}\label{f-Ap est}
\|f-A_{\varrho, r}(f)\|_p={\mathcal O}\left((1-\varrho)^{r-n}\omega(1-\varrho)\right),\quad\varrho\to 1-.
\end{equation}
}

{\theorem\label{Th2}  Assume that $f\in L_{p,Y}(\mathbb T^d),$ $1\le p\le\infty$, $n, r\in\mathbb N$, $n\le r$ and the function  $\omega(t)$, $t\in [0,1]$, satisfies conditions 1)--4), $({\mathcal Z})$  and $({\mathcal Z}_n)$ . If relation (\ref{f-Ap est}) holds, then  relations (\ref{K-funct est})   hold as well.
}

{\remark  For a given $n\in\mathbb N$, from condition $({\mathcal  Z}_n)$  it follows that  $\mathop{\rm lim~inf}\limits_{\delta\to 0+}(\delta^{-n}\omega(\delta))>0$ or, equivalently, that $(1-\varrho)^{r-n}\omega(1-\varrho)\gg (1-\varrho)^r$ as $\varrho\to~1-$. Therefore, if condition $({\mathcal  Z}_n)$  is satisfied, then the quantity on the right-hand side of (\ref{f-Ap est}) decreases to zero as $\varrho{\to} 1-$ not faster than the function $(1-\varrho)^r$. Also note that the relation $
\|f-A_{\varrho, r}(f)\|_p=$\mbox{\tiny  $\mathcal O$}$\left((1-\varrho)^r)\right),\ \varrho\to 1-,
$
 only holds  in the trivial case when $f({\bf x})=\sum_{\nu=0}^{n-1}\sum_{|{\bf k}|_1=\nu}\widehat f_{\bf k}{\mathrm e}_{\bf k}({\bf x})$, and in such case, the theorems are
 easily 
 true. This fact is   related to the so-called  saturation property of the approximation method  generated by the operator $A_{\varrho, r}$. In particular, in [\ref{Savchuk}], it was shown that the operator $A_{\varrho, r}$  generates the linear approximation method of holomorphic functions which is saturated in the space $H_p$ with the saturation order $(1-\varrho)^r$ and the saturation class $H^{r-1}_p\mathop{\rm Lip}1$.}

\subsection{Norms of multipliers in the spaces $L_{p,Y}(\mathbb T^d)$}\label{Norms_mult}
Before proving   Theorems \ref{Th1} and \ref{Th2}, let us give some auxiliary results. In particular,
the following  
Lemma \ref{Lemma_Norm} shows that norms of multipliers in the spaces $L_{p,Y}(\mathbb T^d)$ are equivalent  for all $d$. In our opinion, such a result is interesting in itself.

Let ${\rm M}=\{\mu_\nu\}_{\nu=0}^\infty$ be a sequence of arbitrary complex numbers. If, for any  function $f\in L_{1,Y}(\mathbb T^d)$ with   Fourier series of the form (\ref{Fourier series}), there exists a function $g\in L_{1,Y}(\mathbb T^d)$ with   Fourier series of the form
$$
S[g]({\bf x})=\sum_{\nu=0}^{\infty}\mu_\nu \sum_{{\bf k}\in Y:\,|{\bf k}|_1=\nu}\widehat f_{\bf k} {\mathrm e}_{\bf k}({\bf x}),
$$
then we say that in the space  $L_{1,Y}(\mathbb T^d)$ the multiplier ${\rm M}$ is defined. In this case we use the notation $g={\rm M} (f)$.


Let $B_{p,Y}$, $1\le p\le\infty$, be a unit ball of the space $L_{p,Y}(\mathbb T^d)$, that is, the set of all functions $f\in
L_{p,Y}(\mathbb T^d)$ such that $\|f\|_{_{\scriptstyle p}}\le 1.$

If ${\rm M}: L_{p,Y}(\mathbb T^d) \rightarrow L_{p,Y}(\mathbb T^d)$, 
then the norm of the operator  ${\rm M}$ is the number
$$
\|{\rm M}\|_{_{ _{\scriptstyle L_{p,Y}(\mathbb T^d)\rightarrow L_{p,Y}(\mathbb T^d)}}}=\sup_{f\in
B_{p,Y}}\|{\rm M}(f)\|_{_{\scriptstyle p}}
=\sup_{f\in L_{p,Y}(\mathbb T^d),\atop f\neq 0} \frac{\|{\rm M}(f)\|_{_{\scriptstyle p}}}{\|f\|_{_{\scriptstyle p}}}.
$$
We also denote by $\|{\rm M}\|_{_{ _{\scriptstyle L_p(\mathbb T^1)\rightarrow
L_p(\mathbb T^1)}}}$ the norm of the  operator ${\rm M} :L_p(\mathbb T^1)\rightarrow L_p(\mathbb
T^1)$.

Let us note that if ${\rm M}$ is a continues  operator from
$L_{p,Y}(\mathbb T^d)$ to $L_{p,Y}(\mathbb T^d)$, then
${\rm M}$  is called the multiplier of series  of the form (\ref{Fourier series}) of
$(p,p)$-type (see, for example, [\ref{Edvards}, Ch.~16]).

In [\ref{SavchukV_SavchukM}], the authors proved that
the norms of the multipliers  ${\rm M}$, which are defined in a similar way,
for the Hardy spaces $H_p(\mathbb D^d)$ and $H_p(\mathbb D^1)$ are equivalent   for all $d\in\mathbb N$. Without going into the details,
 we note that the space  $H_p(\mathbb D^d)$ can be considered as the space of all complex-valued
  functions  $f:{\mathbb T^d}\rightarrow\mathbb C$ such that $|f|\in L_p({\mathbb T^d})$ and
$\widehat f({\bf k})=0$ for all  ${\bf k}\in\mathbb Z^d\setminus\mathbb Z^d_{+}$ (see, for example, Theorem 2.1.4
[\ref{Rudin}]). Here, we complement the result of [\ref{SavchukV_SavchukM}] and show that the norms of the multipliers ${\rm M} : L_{p,Y}(\mathbb T^d)\rightarrow
L_{p,Y}(\mathbb T^d)$ are equal as well.


{\lemma\label{Lemma_Norm} Assume that  $1\le p\le\infty$, $d\in\mathbb N$ and ${\rm M}$ is a multiplier
generated by a sequence of complex numbers  $\{\mu_\nu\}_{\nu=0}^\infty$. Then
 \begin{equation}\label{Lemma_norm_equal}
\|{\rm M}\|_{_{_{\scriptstyle L_{p,Y}(\mathbb T^d)\rightarrow L_{p,Y}(\mathbb T^d)}}}=\|{\rm M}\|_{_{
_{\scriptstyle L_{p}(\mathbb T^1)\rightarrow L_{p}(\mathbb
T^1)}}}.
\end{equation}}

\noindent{\it Proof.} Let $f\in L_{p,Y}(\mathbb T^d)$. Note that  for almost all ${\bf x}\in{\mathbb T^d}$, the multiplier ${\rm M}$ can be defined by the  following rule:
 \begin{equation}\label{DEFFFF}
{\rm M}(f)({\bf x})= \lim_{\varrho\to 1-}{\rm M}(f)(\varrho,{\bf x}),
\end{equation}
 where for $0<\varrho<1$ and  ${\bf x}\in{\mathbb T^d}$,
$$
{\rm M}(f)(\varrho,{\bf x})=\sum_{\nu=0}^{\infty}\lambda_\nu\varrho^{\nu}\sum_{{\bf k}\in Y:\,|{\bf k}|_1=\nu}\widehat f_{\bf k} {\mathrm e}_{\bf k}({\bf x}).
$$
If  $f\in L_p(\mathbb T)$, then this rule has the form
$$
{\rm M}(f)(\varrho,t)=\lim_{\varrho\to 1-}\sum_{n\in\mathbb
Z}\mu_{|n|}\varrho^{|n|}\widehat
f_{n}{\mathrm e}^{{\mathrm i}nt}.
$$
For any $f\in L_{p,Y}(\mathbb T^d)$, we set ${\rm M}(f)({\bf z})={\rm M}(f)(\mbox{\boldmath${\bf \varrho}$},{\bf x})$, where
for $0<\varrho_j<1$ and  ${\bf x}\in{\mathbb T^d}$, the point ${\bf z}:=(\varrho_1 {\rm e}^{{\mathrm i}x_1},\ldots,
\varrho_d {\rm e}^{{\mathrm i}x_d})$ belongs to the
unit polydisc $\mathbb D^d:=\{{\bf z}\in{\mathbb C}^d : \max_{1\le j\le d}|z_j|<1\}$. Therefore, the function  ${\rm M}(f)({\bf z})$ is a $d$--harmonic function in ${\mathbb D}^d$ and according to the assertion (c) of Theorem 2.1.3 [\ref{Rudin}], we have
$\|{\rm M}(f)(\varrho{\bf \cdot})\|_{_{\scriptstyle p}}\le\|{\rm M}(f) \|_{_{\scriptstyle p}}$. On the other hand, by virtue of Fatou's lemma,
\[
\|{\rm M}(f) \|_{_{\scriptstyle p}}\le\mathop{\rm lim~inf}\limits_{\varrho\to 1-}\|{\rm M}(f)(\varrho{\bf
\cdot})\|_{_{\scriptstyle p}},
\]
 hence, for  $1\le p<\infty$,
\begin{equation}\label{p4}
\|{\rm M}(f) \|_{_{\scriptstyle p}}=\lim_{\varrho\to 1-}\|{\rm M}(f)(\varrho,{\bf \cdot})\|_{_{\scriptstyle p}}.
\end{equation}
If $p=\infty$, then  instead of the last relation we have
\[
\int\limits_{{\mathbb T^d}}{\rm M}(f)({\bf w})g({\bf w}){\mathrm d}\sigma ({\bf w})
=\lim_{\varrho\to 1-}\int\limits_{{\mathbb T^d}}{\rm M}(f)(\varrho,{\bf w})g({\bf
w}){\mathrm d}\sigma ({\bf w})
\]
for any function  $g\in L_1({\mathbb T^d})$, i.e., we have convergence in the weak $L_1$--topology of space
$L_{\infty}({\mathbb T^d})$.


Let  $f\in L_{p,Y}(\mathbb T^d)$, $f\not\equiv 0$, 
 ${\bf z}$ be a fixed point in $\bar{\mathbb D}^d$ and $0\le \varrho<1$. In the disc $\mathbb D^1$,
consider the function  $u_{\varrho{\bf z}}(\omega):=f(\varrho,{\bf z}\omega)$.
Applying Lemma 3.3.2  [\ref{Rudin}], we consistently have the following equality and estimate for the integral of
$|{\rm M} (f)(\varrho\,\cdot)|^p$ for  $0\le \varrho<1$ and $1\le p<\infty$:
\begin{eqnarray}\nonumber
\lefteqn{\int\limits_{{\mathbb T^d}}|{\rm M} (f)(\varrho, {\bf w})|^p{\mathrm d}\sigma ({\bf w})=
\int\limits_{{\mathbb T^d}}{\mathrm d}\sigma ({\bf w})\int\limits_{\mathbb
T^1}|{\rm M} (u_{\varrho {\bf w}})(\omega)|^pd \omega} \\ \nonumber
&=&\int\limits_{{\mathbb T^d}}\|{\rm M} (u_{\varrho {\bf w}})\|_{_{\scriptstyle p}}^p{\mathrm d}\sigma ({\bf w})
=\int\limits_{{\mathbb T^d}}\|u_{\varrho {\bf w}}\|^p_{_{\scriptstyle p}}\frac{\|{\rm M} (u_{\varrho {\bf w}})\|_{_{\scriptstyle p}}^p}{\|u_{\varrho {\bf w}}\|^p_{_{\scriptstyle p}}}{\mathrm d}\sigma ({\bf w})\\ \nonumber
&\le&\max_{{\bf w}\in \mathbb T^d}\frac{\|{\rm M} (u_{\varrho {\bf w}})\|_{_{\scriptstyle p}}^p}{\|u_{\varrho {\bf w}}\|^p_{_{\scriptstyle p}}}\int\limits_{{\mathbb T^d}}\|u_{\varrho {\bf w}}\|^p_{_{\scriptstyle p}}{\mathrm d}\sigma ({\bf w})
\\ \nonumber
&\le&\|{\rm M} \|^p_{_{_{\scriptstyle L_{p}(\mathbb
T^1)\rightarrow L_{p}(\mathbb T^1)}}}\int\limits_{{\mathbb T^d}}\|u_{\varrho {\bf w}}
\|^p_{_{\scriptstyle p}}{\mathrm d}\sigma ({\bf w})\\ \label{p<infty}
&=&\|{\rm M} \|^p_{_{_{\scriptstyle L_{p}(\mathbb T^1)\rightarrow
L_{p}(\mathbb T^1)}}}\int\limits_{{\mathbb T^d}}|f(\varrho,  {\bf w})|^p{\mathrm d}\sigma ({\bf w}).
\end{eqnarray}
In the case $p=\infty$, we similarly obtain the estimate
\begin{eqnarray}\nonumber
\lefteqn{|{\rm M}(f)(\varrho,\omega{\bf w})|=|{\rm M} (u_{\varrho {\bf w}})(\omega)|}\\ \nonumber &=&
\lim_{\rho\to 1-}|{\rm M} (u_{\varrho {\bf w}})(\rho\omega)|
\le \max_{\omega\in\mathbb T^1}|{\rm M} (u_{\varrho {\bf w}})(\omega)|
\\ \label{p=infty}
&\le&\|{\rm M} \|_{_{_{\scriptstyle L_{\infty,Y}({\mathbb T^d})\rightarrow L_{\infty,Y}({\mathbb T^d})}}}\max_{\omega\in\mathbb T^1}|f(\varrho,   \omega{\bf w})|.
\end{eqnarray}
From 
(\ref{p<infty}) and (\ref{p=infty}) in view of (\ref{DEFFFF}) it follows that for  $1\le
p\le\infty$,
\begin{eqnarray}\nonumber
\|{\rm M} \|_{_{ _{\scriptstyle L_{p,Y}(\mathbb T^d)\rightarrow L_{p,Y}(\mathbb T^d)}}}&=&
\lim_{\varrho\to 1-}\sup_{f\in L_{p,Y}(\mathbb T^d)}\frac{\|{\rm M} (f)(\varrho,\cdot)
\|_{_{\scriptstyle p}}}{\|f(\varrho,\cdot)\|_{_{\scriptstyle p}}}
\\ \label{upper_est}
&\le&\|{\rm M} \|_{_{_{\scriptstyle L_{p}(\mathbb
T^1)\rightarrow L_{p}(\mathbb T^1)}}}.
\end{eqnarray}

To prove the reverse inequality let us consider the continuation operator
$Q$, given on  $L_{p}(\mathbb T^1)$, $1\le p\le\infty$, by the formula
$$
Q(g)(w_1,{\bf w}^1)=g(w_1),
$$
where $w_1\in\mathbb T^1,~{\bf w}^1=(w_2,\ldots,w_d)\in\mathbb T^{d-1}$.

It is easy to show that the continuation operator $Q$ is a linear isometry of the space $L_{p}(\mathbb T^1)$ in $L_{p}({\mathbb T^d})$. Therefore, taking into account the relation
$Q\big({\rm M} (f)\big)={\rm M} \big(Q(f)\big)$, which is satisfied for any function $f\in L_p(\mathbb T^1)$,
we obtain
\begin{eqnarray}\nonumber
\lefteqn{\|{\rm M} (f)\|_{_{\scriptstyle p}}=\|Q\big({\rm M} (f)\big)\|_{_{\scriptstyle p}}=\|{\rm M} \big(Q(f)\big)\|_{_{\scriptstyle p}}}
 \\ \nonumber
&\le&\|{\rm M} \|_{_{_{\scriptstyle L_{p,Y}(\mathbb T^d)\rightarrow L_{p,Y}(\mathbb T^d)}}}\|Q(f)\|_{_{\scriptstyle p}}=\|{\rm M} \|_{_{_{\scriptstyle L_{p,Y}(\mathbb T^d)\rightarrow L_{p,Y}(\mathbb T^d)}}}\|f\|_{_{\scriptstyle p}}.
\end{eqnarray}
This implies the estimate
$$
\|{\rm M} \|_{_{ _{\scriptstyle L_{p}(\mathbb T^1)\rightarrow
L_{p}(\mathbb T^1)}}}\le\|{\rm M} \|_{_{_{\scriptstyle
L_{p,Y}(\mathbb T^d)\rightarrow L_{p,Y}(\mathbb T^d)}}},
$$
 which in combination with (\ref{upper_est}) gives the relation (\ref{Lemma_norm_equal}).


\subsection{Auxiliary statements}\label{Subsetion2.3}
  Let
 \begin{equation}\label{Poisson operator_NEW}
{\mathcal P}(\varrho,{\bf x}):=\prod_{j=1}^d\frac 1{1-\varrho {\mathrm e}^{{\mathrm i}x_j}}+\prod_{j=1}^d\frac 1{1-\varrho {\mathrm e}^{-ix_j}}-1.
\end{equation}

{\lemma\label{Lemma0} Assume that $f\in L_{1,Y}(\mathbb T^d)$, $0\le\varrho<1$ and ${\bf x}\in \mathbb T^d$. Then
 \begin{equation}\label{Poisson operator2}
f\left(\varrho,{\bf x}\right)=  \int_{{\mathbb T}^d}
f({\bf x}+{\bf s}){\mathcal P}(\varrho,{\bf s}){\mathrm d}\sigma({\bf s}).
\end{equation}
}

\noindent{\it Proof.}  By virtue of the definition of the set  $L_{1,Y}(\mathbb T^d)$, we have
 \begin{equation}\label{Poisson operator_represent}
f(\varrho ,{\bf x})=\sum_{\nu=0}^{\infty}\varrho^{\nu}\sum_{{\bf k}\in Y:\,|{\bf k}|_1=\nu}\widehat f_{\bf k} {\mathrm e}_{\bf k}({\bf x}).
\end{equation}
 On the other hand
 \begin{eqnarray}\nonumber
{\mathcal P}(\varrho,{\bf x})&=&\sum\limits_{k_1=0}^\infty
\ldots\sum\limits_{k_d=0}^\infty
\varrho^{k_1+
\ldots+k_d} \Big({\mathrm e}^{{\mathrm i}(k_1x_1+
\ldots+k_dx_d)}+
{\mathrm e}^{-{\mathrm i}(k_1x_1+
\ldots+k_dx_d)}\Big)-1\\ \label{Poisson operator_represent_new}
&=&
1+\sum\limits_{\nu=1}^\infty\varrho ^{{\nu}}\sum\limits_{{\bf k}\in Y:\,|{\bf k}|_1=\nu}
{\mathrm e}_{\bf k}({\bf x}).
\end{eqnarray}
Therefore, the right-hand side of (\ref{Poisson operator2}) is equivalent to the right-hand side  of  (\ref{Poisson operator_represent}).

 \ $\hfill\Box$


{\lemma\label{Lemma3}  Assume that  $f\in L_{p,Y}(\mathbb T^d)$, $1\le p\le \infty$,  $r=0,1,\ldots$  and $\varrho\in [0,1)$. Then the following relations are true:
\begin{equation}\label{|derivative_Poisson_Integral|}
\Big\| \frac{\partial^r f\left(\varrho,\cdot\right)}{\partial\varrho^r}\Big\|_{p}\le C_{1}(r)\frac{\|f\|_p}{(1- \varrho
)^{r}}
\end{equation}
and
\begin{equation}\label{|A_r(f)|}
\|A_{\varrho, r}^{[r]}(f)\|_p \le C_{2}(r)\frac{\|f\|_p}{(1-\varrho)^r},
\end{equation}
where the
constants $C_{1}(r)$ and $C_{2}(r)$ depend only on $r$.
}

\noindent {\it Proof.} It is easy to see that the function $\partial^rf\left(\varrho,{\bf x}\right)/ {\partial\varrho^r}$  can be considered as the image ${\rm M}_1(f)({\bf x})$  of  the multiplier generated by the sequence $\{\mu_{1,\nu}\}_{\nu=0}^\infty$, where $\mu_{1,\nu}=0$ for $\nu=0,1,\ldots, r-1$ and
 $\mu_{1,\nu}=\nu\cdot(\nu-1)\cdot\ldots\cdot (\nu-r+1)\varrho^{\nu-r}$ for $\nu\ge r$. Similarly,  the function $A_{\varrho, r}^{[r]}(f)({\bf x})$ can be considered as the image ${\rm M}_2(f)({\bf x})$  of the multiplier generated by the sequence $\{\mu_{2,\nu}\}_{\nu=0}^\infty$ such that $\mu_{2,\nu}=0$ for $\nu=0,1,\ldots, r-1$ and
 $\mu_{2,\nu}=\nu! \cdot \lambda_{\nu,r}(\varrho)/ (\nu-r)!$ for $\nu\ge r$. Therefore, to prove estimates
 (\ref{|derivative_Poisson_Integral|}) and (\ref{|A_r(f)|}) it is sufficient to apply Lemma \ref{Lemma_Norm} and
  the estimates (23) and (22) for the norms of
 the corresponding multipliers in the space $L_{p}(\mathbb T^1)$ from [\ref{Prestin_Savchuk_Shidlich}].  \vskip -3mm$\hfill\Box$


For any $f\in L_p(\mathbb T^d)$,  $1\le p \le \infty$, $0\le\varrho<1$ and $r=0,1,2,\ldots$, we set
\begin{equation}\label{M_p_diff}
M_p(\varrho, f,r):=\varrho^r\Big\| \frac{\partial^r f\left(\varrho,\cdot\right)}{\partial\varrho^r}\Big\|_{p}=
\Big\|(f(\varrho,\cdot))^{[r]} \Big\|_{p}.
\end{equation}

{\lemma\label{Main_Lemma} Assume that  $f\in L_{p,Y}(\mathbb T^d)$, $1\le p\le \infty$. Then for any numbers $n\in \mathbb N$ and $\varrho\in[0,1)$,
\begin{eqnarray}\nonumber
\lefteqn{C_{3}(n)(1- \varrho)^{n}M_p\left(\varrho, f,n\right)\le K_{n}\left(1- \varrho, f\right)_p}
\\ \label{Ineq_Lemma}
&\le&  C_{4}(n)\Big(\|f-A_{\varrho, n}(f)\|_p+(1-\varrho)^{n}M_p\left(\varrho, f,n\right)\Big),
\end{eqnarray}
where  the constants $C_{3}(n)$ and $C_{4}(n)$  depend only on $n$.
}

\noindent {\it Proof.}   First, let us note that the statement of Lemma \ref{Main_Lemma} is trivial in the case,
when
$f$ is a polynomial of the form $f({\bf x})=\sum_{\nu=0}^{n-1}\sum_{|{\bf k}|_1=\nu}\widehat f_{\bf k}{\mathrm e}_{\bf k}({\bf x})$, as well as in the case,
when
 $\varrho=0$. Therefore, further in the proof, we exclude these two cases.

Let $g$ be a function such that  $g^{[n]}\in L_p(\mathbb T^d)$. Using Lemma \ref{Lemma3}, we get
\begin{eqnarray*}
\Big\|\frac{\partial^n f\left(\varrho,\cdot\right)}{\partial\varrho^n}\Big\|_p&=&
\Big\|\frac{\partial^n (f-g)\left(\varrho,\cdot\right)}{\partial\varrho^n}+ \frac{\partial^n g\left(\varrho,\cdot\right)}{\partial\varrho^n}\Big\|_p\\ &\le&  C_1(n)\frac{\|f-g\|_p}
{(1- \varrho)^{n}}+\Big\|\frac{\partial^n g\left(\varrho,\cdot\right)}{\partial\varrho^n}\Big\|_p.
\end{eqnarray*}
Setting $C_{3}(n)=\min\{1,1/C_1(n)\}$ and taking into account relations (\ref{diff f[n]}),  (\ref{M_p_diff})  and the inequality $\|g^{[n]}(\varrho,\cdot)\|_p\le\|g^{[n]}\|_p$, we see that
 \[
C_{3}(n)(1- \varrho)^{n}M_p(\varrho, f,n)\le \|f-g\|_p+(1- \varrho)^{n} \|g^{[n]}\|_p.
\]
Considering the infimum over all functions $g$ such that $g^{[n]}\in L_p(\mathbb T^d)$, we conclude that
\[
C_{3}(n)(1-\varrho)^nM_p\left(\varrho, f,n\right)\le K_{n}\left(1- \varrho, f\right)_p.
\]
On the other hand, from the definition of the  $K$--functional, it follows that
\begin{equation}\label{K_n}
K_n\left(1- \varrho, f\right)_p\le\|f-A_{\varrho, n}(f)\|_p+(1-\varrho)^{n}\left\|
A_{\varrho,n}^{[n]}(f)\right\|_p.
\end{equation}
According to (\ref{A_P}) and (\ref{diff f[n]}), we have
\begin{eqnarray*}
A_{\varrho,n}^{[n]}(f)({\bf x})&=&\bigg(\sum_{k=0}^{n-1}\frac{(f(\varrho,\cdot))^{[k]}(\cdot)}{\varrho^k k!}(1-\varrho)^k\bigg)^{[n]}({\bf x})
\\ &=&\sum_{k=0}^{n-1}\frac{((f(\varrho,\cdot))^{[k]}(\cdot))^{[n]}({\bf x})}{\varrho^k k!}(1-\varrho)^k.
\end{eqnarray*}

Since for any nonnegative integers $k$ and $n$
\begin{equation}\label{diff f[n]_f[k]}
((f(\varrho,\cdot))^{[n]}(\cdot))^{[k]}({\bf x})=((f(\varrho,\cdot))^{[k]}(\cdot))^{[n]}({\bf x}),
\end{equation}
we obtain
$$
A_{\varrho,n}^{[n]}(f)({\bf x})=\sum_{k=0}^{n-1}\frac{((f(\varrho,\cdot))^{[n]}(\cdot))^{[k]}({\bf x})}{\varrho^k k!}(1-\varrho)^k.
$$
This yields
\begin{equation}\label{A11}
\|A_{\varrho,n}^{[n]}(f)\|_{_{\scriptstyle p}}\le \sum_{k=0}^{n-1}\frac{\|((f(\varrho,\cdot))^{[n]}(\cdot))^{[k]}\|_{_{\scriptstyle p}}}{\varrho^k k!}(1-\varrho)^k,
\end{equation}
where by virtue of Lemma \ref{Lemma3} and (\ref{M_p_diff})
\begin{equation}\label{Ab1}
\|((f(\varrho,\cdot))^{[n]}(\cdot))^{[k]}\|_{_{\scriptstyle p}}
\le
M_p(\varrho,f,n)\frac{C_{1}(k)\varrho^k}{(1-\varrho)^{k}}.
\end{equation}
Therefore,
\begin{equation}\label{Ab2}
\|A_{\varrho,n}^{[n]}(f)\|_{_{\scriptstyle p}}\le M_p(\varrho,f,n)\sum_{k=0}^{n-1}\frac{C_{1}(k)}{k!}.
\end{equation}
Setting $C_{4}(n)=\max\{1,\sum_{k=0}^{n-1}{C_{1}(k)}/{k!}\}$ and combining relations (\ref{K_n}) and (\ref{Ab2}), we obtain the right-hand inequality in (\ref{Ineq_Lemma}). \vskip -3mm$\hfill\Box$



 {\lemma\label{Lemma4}  Assume that $f\in L_p
 (\mathbb T^d)$, $1\le p\le \infty$, $0\le \varrho<1$ and $r=2,3,\ldots$ such that
\begin{equation}\label{f-A(f)Condition}
\int_\varrho^1 \bigg\|\frac{\partial^rf(\zeta,\cdot)}{\partial\zeta^r}\bigg\|_{_{\scriptstyle p}}(1-\zeta)^{r-1}d\zeta<\infty.
\end{equation}
Then for almost all ${\bf x}\in \mathbb T^d$,
\begin{equation}\label{f-A(f)}
f({\bf x})-A_{\varrho, r}(f)({\bf x})=\frac{1}{(r-1)!}\int_\varrho^1 \frac{\partial^rf(\zeta,{\bf x})}{\partial\zeta^r}(1-\zeta)^{r-1}d\zeta.
\end{equation}

}

\noindent {\it Proof.} For  fixed $r=2,3,\ldots$ and $0\le \varrho<1$, the integral on the right-hand side of (\ref{f-A(f)}) defines a certain function $F({\bf x})$. By virtue of (\ref{f-A(f)Condition}) and  the integral Minkowski inequality, we conclude that the function  $F$ belongs to the space $L_p(\mathbb T^d)$. Let us find the Fourier coefficients of $F$ and compare them with  the Fourier coefficients of the function $G:=f-A_{\varrho, r}(f)$.

Since for any $\nu=r,r+1\ldots$,
\[
\frac 1{(r-1)!\cdot(\nu-r)!}\int_\varrho^{\varrho_1}\zeta^{\nu-r}(1-\zeta)^{r-1}d\zeta
=\sum_{j=0}^{r-1}
\frac{\varrho_1^{\nu-j}(1-\varrho_1)^{j}-
\varrho^{\nu-j}(1-\varrho)^{j}}{j!\cdot(\nu-j)!},
\]
then in view of  (\ref{Poisson_Derivative})  for a fixed $\varrho_1\in (\varrho,1)$, we have
\begin{eqnarray}\nonumber
\lefteqn{ \frac 1{(r-1)!}\int_\varrho^{\varrho_1} \frac{\partial^rf(\zeta,{\bf x})}{\partial\zeta^r}(1-\zeta)^{r-1}{\mathrm d}\zeta}\\ \nonumber
&=&
\sum_{\nu=r}^{\infty}\sum_{|{\bf k}|_1=\nu}\frac{ \nu!\cdot \widehat f_{\bf k}\cdot {\mathrm e}_{\bf k}({\bf x})}{(r-1)!\cdot (\nu-r)!}
\int_\varrho^{\varrho_1} \zeta^{\nu-r}(1-\zeta)^{r-1}{\mathrm d}\zeta
\\ \label{Ab3}
&=&\sum_{\nu=r}^{\infty}\sum_{|{\bf k}|_1=\nu} \widehat f_{\bf k} {\mathrm e}_{\bf k}({\bf x})
\sum_{j=0}^{r-1}
{\nu\choose j} \Big(\varrho_1^{\nu-j}(1-\varrho_1)^{j}-
\varrho^{\nu-j}(1-\varrho)^{j}\Big).
\end{eqnarray}
Now if in relation (\ref{Ab3}), the value $\varrho_1$ tends to $1-$, then  we see that the Fourier coefficients $\widehat{F}_{{\bf k}}$  of the  function $F$ are equivalent to zero when $|{\bf k}|_1=\nu<r$ and for $|{\bf k}|_1\ge r$,
\begin{equation}\label{Fourier coeff_F}
\widehat{F}_{{\bf k}}=\widehat f_{\bf k}\cdot \Big(1-\sum_{j=0}^{r-1}
{\nu\choose j} (1-\varrho)^j\varrho^{\nu-j}\Big)=(1-\lambda_{\nu,r}(\varrho))\widehat f_{\bf k}.
\end{equation}
Therefore, for all ${\bf k}\in {\mathbb Z}^d$ we have  $\widehat{F}_{{\bf k}}=(1-\lambda_{\nu,r}(\varrho))\widehat f_{\bf k}=\widehat{G}_{\bf k}$. Hence, for almost all ${\bf x}\in \mathbb T^d$, relation (\ref{f-A(f)}) holds.  \vskip -3mm$\hfill\Box$

\subsection{Proof of main results}
  {\it  Proof of Theorem \ref{Th1}.} Assume that the function $f$ is such that $f^{[r-n]}\in L_{p,Y}(\mathbb T^d)$ and relation (\ref{K-funct est}) is satisfied.  Let us apply  the first inequality
of
Lemma \ref{Main_Lemma} to the function $f^{[r-n]}$. In view of  (\ref{diff f[n]}) and (\ref{M_p_diff}), we obtain
\[
C_{3}(n)(1-\varrho)^{n}M_p\left(\varrho, f,r\right)\le K_{n}\left(1-\varrho, f^{[r-n]}\right)_p.
\]
This yields
\begin{equation}\label{est Mp}
M_p\left(\varrho, f,r\right)={\mathcal O}(1) {(1-\varrho)^{-n}}{\omega(1-\varrho)}, \quad \varrho\to 1-.
\end{equation}

Using  relations (\ref{M_p_diff}), (\ref{est Mp}), $({\mathcal Z})$  and the integral Minkowski inequality, we obtain
\begin{equation}\label{est Mp estim}
 \int_\varrho^1 \bigg\|\frac{\partial^rf(\zeta,\cdot)}{\partial\zeta^r}\bigg\|_{_{\scriptstyle p}}(1-\zeta)^{r-1}{\mathrm d}\zeta \le \int_\varrho^1 M_p\left(\zeta, f,r\right)\frac{(1-\zeta)^{r-1}}{\zeta^r}{\mathrm d}\zeta
 $$
 $$
 \le
 {C_1}(1-\varrho)^{r-n}\int_\varrho^1\frac{\omega(1-\zeta)}{1-\zeta}{\mathrm d}\zeta=
 {\mathcal O}\left((1-\varrho)^{r-n}\omega(1-\varrho)\right), \ \varrho\to 1-.
\end{equation}
Therefore, for almost all ${\bf x}\in \mathbb T^d$, relation (\ref{f-A(f)}) holds. Hence, by virtue of  (\ref{f-A(f)}),
using the integral Minkowski inequality and (\ref{est Mp estim}),  we
finally 
get (\ref{f-Ap est}):
 \begin{eqnarray*}
\|f-A_{\varrho, r}(f)\|_{_{\scriptstyle p}}&\le&\frac{1}{(r-1)!}\int_\varrho^1 M_p\left(\zeta, f,r\right)\frac{(1-\zeta)^{r-1}}{\zeta^r}{\mathrm d}\zeta \\ &=&{\mathcal O}\left((1-\varrho)^{r-n}\omega(1-\varrho)\right), \quad \varrho\to 1-.
\end{eqnarray*}
\vskip -3mm$\hfill\Box$

\bigskip
\noindent {\it  Proof of Theorem \ref{Th2}.} First, let us note that for any function $f\in L_p(\mathbb T^d)$ and all fixed numbers $s,r\in {\mathbb N}$ and $\varrho\in (0,1)$, we have

\begin{eqnarray}\nonumber
\|A_{\varrho, r}^{[s]}(f)\|_{_{\scriptstyle p}}&=&\Big\|
 \sum_{\nu=s}^{\infty}\frac{\nu!}{(\nu-s)!}\lambda_{\nu,r}(\varrho)\sum_{|{\bf k}|_1=\nu}\widehat f_{\bf k}{\mathrm e}_{\bf k} \Big\|_{_{\scriptstyle p}}
 \\ \nonumber
&\ \le&   2r\|f\|_{_{\scriptstyle p}}\bigg(\sum_{\nu=s}^{\max\{s,r\}-1}\frac{\nu!}{(\nu-s)!}+\sum_{\nu\ge \max\{s,r\}}q^{\nu}\nu^{s+r-1}\bigg)<\infty,
\end{eqnarray}
where $0<q=\max\{1-\varrho,\varrho\}<1$. In the case where $s\ge r$,  the  sum $\sum_{\nu=s}^{s-1}$ is set equal to zero.

Put $\varrho_k:=1-2^{-k},~k\in{\mathbb N},$ and $A_k:=A_k(f):=A_{\varrho_k,r}(f)$. For any ${\bf x}\in
{\mathbb T}^d$
and $s\in {\mathbb N}$, consider the series
\begin{equation}\label{series1}
A_0^{[s]}(f)({\bf x})+\sum\limits_{k=1}^\infty (A_k^{[s]}(f)({\bf x})-A_{k-1}^{[s]}(f)({\bf x})).
\end{equation}
According to the definition of the operator  $A_{\varrho,r}$, we see that for any $\varrho_1, \varrho_2\in[0,1)$ and $r\in {\mathbb N}$,
\[
A_{\varrho_1,r}\left(A_{\varrho_2,r}(f)\right)=A_{\varrho_2,r}\left(A_{\varrho_1,r}(f)\right).
\]
By virtue of Lemma \ref{Lemma3} and relation (\ref{f-Ap est}), for any  $k\in {\mathbb N}$ and  $s\in {\mathbb N}$, we have
\begin{eqnarray}\nonumber
\lefteqn{\left\|A^{[s]}_{k}-A^{[s]}_{k-1}\right\|_{_{\scriptstyle p}}
=\left\|A^{[s]}_{k}(f-A_{k-1}(f))-A^{[s]}_{k-1}(f-A_{k}(f))
\right\|_{_{\scriptstyle p}}}\\ \nonumber
&\le&\left\|A^{[s]}_{k}(f-A_{k-1}(f))\right\|_{_{\scriptstyle p}}+\left\|A^{[s]}_{k-1}(f-A_{k}(f)) \right\|_{_{\scriptstyle p}}\\ \nonumber
&\le& C_{2}(s)\frac{\left\|f-A_{k-1}(f)\right\|_{_{\scriptstyle p}}}{(1-\varrho_k)^{s}}+C_{2}(s)\frac{\left\|f-A_{k}(f)\right\|_{_{\scriptstyle p}}}
{(1-\varrho_{k-1})^{s}} \\ \label{OA}
&=&
{\mathcal O}\left(\frac{\omega(1-\varrho_{k-1})}{(1-\varrho_{k})^{s-r+n}}\right)+
{\mathcal O}\left(\frac{\omega(1-\varrho_{k})}{(1-\varrho_{k-1})^{s-r+n}}\right),\quad k\to \infty.
\end{eqnarray}
Therefore, for any $s\le r-n$,
\begin{equation}\label{series12}
\left\|A^{[s]}_{k}-A^{[s]}_{k-1}\right\|_{_{\scriptstyle p}}={\mathcal O}\left(\omega(1-\varrho_{k-1})\right)=
{\mathcal O}\left(\omega(2^{-(k-1)})\right),\quad k\to \infty.
\end{equation}
Consider the sum  $\sum_{k=1}^N \omega(2^{1-k})$, $N\in {\mathbb N}$. Taking into account the monotonicity of the function $\omega$ and $({\mathcal Z})$ , we see that for all $N\in {\mathbb N}$,
\[
 \sum_{k=1}^N\omega(2^{1-k})\le \omega(1)+\int_1^{N} \omega(2^{1-t}){\mathrm d}t=
 \omega(1)+ \int_{2^{1-N}}^1 \frac{\omega(\tau)\,d\tau}{\tau\ln 2} <\infty.
\]

Combining the last relation and (\ref{series12}), we  conclude that  the series in (\ref{series1}) converges in the
norm
of the space $L_p(\mathbb T^d)$, $1\le p\le\infty$. Hence, by virtue of the Banach--Alaoglu theorem, for any  $s=0,1,\ldots,  r-n$, there exists the subsequence
\begin{equation}\label{series123}
S^{[s]}_{N_j}({\bf x})=A_0^{[s]}(f)({\bf x})+\sum\limits_{k=1}^{N_j} (A_k^{[s]}(f)({\bf x})-A_{k-1}^{[s]}(f)({\bf x})),\quad j=1,2,\ldots
\end{equation}
of partial sums of this series, converging to a certain function  $g\in L_p(\mathbb T^d)$ almost everywhere on
${\mathbb T}^d$
 as $j\to\infty$.

Let us show that $g=f^{[s]}$. For this, let us find the Fourier coefficients of the function $g$. For any fixed  ${\bf k}\in{\mathbb Z}^d$ and all $j=1,2,\ldots,$ we have
\[
\widehat g_{\bf k}:=\int_{{\mathbb T}^d}S^{[s]}_{N_j}({\bf x})\overline{\mathrm e}_{\bf k}({\bf x}){\mathrm d}\sigma({\bf x})+
\int_{{\mathbb T}^d}(g({\bf x})-S^{[s]}_{N_j}({\bf x}))\overline{\mathrm e}_{\bf k}({\bf x}){\mathrm d}\sigma({\bf x}).
\]
Since the sequence $\{S^{[s]}_{N_j}\}_{j=1}^\infty$ converges almost everywhere on
${\mathbb T}^d$
to the function $g$, 
the second integral on the right-hand side of the last equality tends to zero as $j\to\infty$. By virtue of (\ref{series123}) and the definition of the radial derivative,
for $|{\bf k}|_1=\nu<s$ the first integral is equal to zero,
and for all $|{\bf k}|_1=\nu\ge s$,
 $$
\int_{{\mathbb T}^d}S^{[s]}_{N_j}({\bf x})\overline{\mathrm e}_{\bf k}({\bf x}){\mathrm d}\sigma({\bf x})=
 \lambda_{\nu,r}(1-2^{-N_j})\frac{\nu!}{(\nu-s)!}\widehat f_{\bf k} \mathop{\longrightarrow}\limits_{j\to\infty} \frac{\nu!}{(\nu-s)!}\widehat f_{\bf k}.
 $$
Therefore, the equality  $g=f^{[s]}$ is true. Hence, for the function  $f$ and all $s=0,1,\ldots,r-n$, there exists the derivative $f^{[s]}$ and $f^{[s]}\in L_p(\mathbb T^d)$.

Now,
let us prove the estimate  (\ref{est Mp}). By virtue of (\ref{M_p_diff}), (\ref{OA}), for any $k\in\mathbb N$ and $\varrho\in (0,1)$, we have
\begin{eqnarray}\nonumber
\lefteqn{M_p\left(\varrho, A_{k}-A_{k-1},r\right)\le \left\|A^{[r]}_{k}-A^{[r]}_{k-1}\right\|_{_{\scriptstyle p}}}
\\ \nonumber
&=&{\mathcal O}\left(\frac{\omega(1-\varrho_{k-1})}{(1-\varrho_{k})^n}\right)+{\mathcal O}\left(\frac{\omega(1-\varrho_{k})}{(1-\varrho_{k-1})^n}\right)
\\ \nonumber
&=&{\mathcal O}\left(2^{kn}\omega(2^{-k+1})+2^{(k-1)n}\omega(2^{-k})\right)
\\ \label{deviation for AA}
&=&
{\mathcal O}\left(2^{(k-1)n}\omega(2^{-(k-1)})\right),\ \  k\to \infty.
\end{eqnarray}
By virtue of  (\ref{M_p_diff}), (\ref{|derivative_Poisson_Integral|})  and (\ref{f-Ap est}),  for any $r\in {\mathbb N}$ and $\varrho\in (0,1)$, we obtain
\[
M_p\left(\varrho, f-A_{\varrho,r}(f),r\right)={\mathcal O}(1) \frac{\left\|f-A_{\varrho,r}(f)\right\|_{p}}{(1-\varrho)^{r}}
={\mathcal O}\left(\frac{\omega(1-\varrho)}{(1-\varrho)^{n}}\right), \quad \varrho\to 1-.
\]
Therefore, for  $N\to \infty$,
\begin{equation}\label{deviation for f}
M_p\left(\varrho_{_{\scriptstyle N}}, f-A_{N}(f),r\right)={\mathcal O}\left(\frac{\omega(1-\varrho_{_{\scriptstyle N}})}{(1-\varrho_{_{\scriptstyle N}})^n}\right)= {\mathcal O}\left(2^{Nn}\omega(2^{-N})\right).
\end{equation}

Consider the sum $\sum_{k=1}^N2^{(k-1)n}\omega(2^{-(k-1)})$, $N\in {\mathbb N}$.  Since the function  $\omega$ satisfies the condition $({\mathcal Z}_n)$ ,
the function $\omega(t)/t^n$ almost decreases on $(0,1]$, i.e., there exists the number $C>0$ such that  $\omega(t_1)/t_1^n\ge C \omega(t_2)/t_2^n$ for any $0<t_1<t_2\le 1$  (see, for example [\ref{Bari_Stechkin}]). Therefore,
\begin{eqnarray}\nonumber
\lefteqn{ \sum_{k=1}^N2^{(k-1)n}\omega(2^{-(k-1)})}\\ \nonumber
&\le& C\bigg( 2^{(N-1)n}\omega(2^{-(N-1)})
 +\int_1^{N} 2^{(t-1)n}\omega(2^{-(t-1)}){\mathrm d}t\bigg)\\ \nonumber
&\le& C\bigg(2^{(N-1)n}\omega(2^{-(N-1)})+\int_{2^{-N+1}}^1\frac{ \omega(\tau)\,d\tau}{
  \tau^{n+1}\ln 2}\bigg)\\ \label{summe}
 &=&{\mathcal O}\Big( 2^{(N-1)n}\omega(2^{-(N-1)})\Big)={\mathcal O}\Big(2^{Nn}\omega(2^{-N})\Big),\quad N\to \infty.
\end{eqnarray}
Putting $\varrho=\varrho_{_{\scriptstyle N}}$ and taking into account relations (\ref{deviation for AA}), (\ref{deviation for f}), (\ref{summe}) and
$$
A_{0}({\bf x})=S_{r-1}(f)({\bf x})=\sum_{|{\bf k}|_1\le r-1}\widehat f_{\bf k}{\mathrm e}_{\bf k}({\bf x}),
$$
we get
\begin{eqnarray}\nonumber
\lefteqn{M_p\left(\varrho_N,f,r\right)=M_p\left(\varrho_N,f-S_{r-1}(f),r\right)}
\\ \nonumber
&=& M_p\bigg(\varrho_N,f-A_{\varrho_{_{\scriptstyle N}}}+\sum\limits_{k=1}^N (A_{k}-A_{k-1}),r\bigg)
= {\mathcal O}\bigg(\sum_{k=1}^{N}2^{(k-1)n}\omega(2^{-(k-1)})\bigg)
\\ \label{M_p(varrho_N)}
&=&  {\mathcal O}\left(2^{Nn}\omega(2^{-N})\right)
 ={\mathcal O}\left((1-\varrho_{_{\scriptstyle N}})^{-n}\omega(1-\varrho_{_{\scriptstyle N}})\right),\quad  N\to \infty.
\end{eqnarray}

If the  function  $\omega$ satisfies the condition $({\mathcal Z}_n)$, then
$\sup_{t\in [0,1]}
\Big({\omega(2t)}/{\omega(t)}\Big)<\infty$  (see, for example [\ref{Bari_Stechkin}]).  Furthermore, for all $\varrho\in[\varrho_{_{\scriptstyle N-1}},\varrho_{_{\scriptstyle
N}}]$, we have $1-\varrho_{_{\scriptstyle N}}\le
1-\varrho\le2(1-\varrho_{_{\scriptstyle N}})$. Hence,  relation (\ref{M_p(varrho_N)}) yields the estimate (\ref{est Mp}).

Now, applying the second inequality in Lemma \ref{Main_Lemma} to the function   $f^{[r-n]}$, we get
\begin{eqnarray}\nonumber
K_{n}\left(1-\varrho, f^{[r-n]}\right)_p&\le& C_{4}(n)\Big( \|f^{[r-n]}-A_{\varrho, n}(f^{[r-n]})\|_{_{\scriptstyle p}}
\\ \label{est Kn}
&+&
(1-\varrho)^nM_p({\varrho},f,r)\Big).
\end{eqnarray}
By virtue of (\ref{M_p_diff}) and (\ref{est Mp}), we see  that for  $\varrho\in [1/2,1)$,
\begin{eqnarray}\nonumber
\lefteqn{\int_\varrho^1 \bigg\|\frac{\partial^n f^{[r-n]}(\zeta,\cdot)}{\partial\zeta^n}\bigg\|_{_{\scriptstyle p}}(1-\zeta)^{n-1}{\mathrm d}\zeta}\\ \nonumber
&=&\int_\varrho^1 \Big\|(f(\zeta,\cdot))^{[r]} \Big\|_{_{\scriptstyle p}}\frac{(1-\zeta)^{n-1}}{\zeta^n} {\mathrm d}\zeta\\ \nonumber
&=&\int_\varrho^1 M_p\left(\zeta, f,r\right)\frac{(1-\zeta)^{n-1}}{\zeta^n} {\mathrm d}\zeta
\\ \label{111} &\le& C_1
\int_\varrho^1\frac{\omega(1-\zeta)}{1-\zeta}{\mathrm d}\zeta
={\mathcal O}\left(\omega(1-\varrho)\right), \quad \varrho\to 1-.
\end{eqnarray}
Therefore, we can apply  Lemma \ref{Lemma4} to the function   $f^{[r-n]}$. Taking into account (\ref{M_p_diff}), we obtain
\[
f^{[r-n]}({\bf x})-A_{\varrho, n}(f^{[r-n]})({\bf x})=\frac{1}{(n-1)!}\int_\varrho^1 (f(\zeta,\cdot))^{[r]}({\bf x})\frac{(1-\zeta)^{n-1}}{\zeta^n} {\mathrm d}\zeta.
\]
Using   the integral Minkowski inequality and (\ref{111}), we 
conclude
\begin{eqnarray}\nonumber
\|f^{[r-n]}-A_{\varrho, n}(f^{[r-n]})\|_{_{\scriptstyle p}}&\le& \frac{1}{(n-1)!}\int_\varrho^1 M_p\left(\zeta, f,r\right)\frac{(1-\zeta)^{n-1}}{\zeta^n} {\mathrm d}\zeta\\ \label{112}
&=& {\mathcal O}\left(\omega(1-\varrho)\right), \quad \varrho\to 1-.
\end{eqnarray}
Combining relations (\ref{est Kn}), (\ref{est Mp}) and (\ref{112}),  we finally get (\ref{K-funct est}).


\vskip 3mm

{\bf\textsl Acknowledgments.} This work is partially supported by the Grant H2020-MSCA-RISE-2014,
project number 645672 (AMMODIT: Approximation Methods for Molecular
Modelling and Diagnosis Tools) and the Grant of the NAS of Ukraine to research laboratories / groups of young scientists of the NAS of Ukraine for conducting researches on priority directions of science and technology development in 2018, project number 16-10/2018.

\vskip 3mm

\end{document}